\documentclass[leqno,10pt]{article}[1999/09/05]

\usepackage{a4,amsfonts}

\newtheorem{satz}{Proposition}

\newtheorem{theorem}[satz]{Theorem}

\newcommand{\F}[1]{ \ensuremath{ \mathbb{F}_{2^{#1}}}}
\newcommand{\Fm}[1]{ \ensuremath{ \mathbb{F}_{2^{#1}}^{\, *}}}
\newcommand{\Fv}[1]{ \ensuremath{ {\mathbb{F}_2^{\; #1}}}}
\newcommand{\C}{ \ensuremath{ \mathbb C}}

\newcommand{\K}{ \ensuremath{ \mathbb K}}

\begin{document}
\bibliographystyle{plain}

\title{ A new APN function which is not equivalent to a power mapping}
\author{
 Yves Edel\footnote{y.edel@mathi.uni-heidelberg.de} \\ Mathematisches Institut der Universit\"at\\
Im Neuenheimer Feld 288\\ D-69120 Heidelberg\\ Germany\\
and\\
Gohar
Kyureghyan\footnote{gohar.kyureghyan@mathematik.uni-magdeburg.de} 
\ and 
Alexander Pott\footnote{alexander.pott@mathematik.uni-magdeburg.de}\\ Institute
for Algebra and Geometry\\ Otto-von-Guericke-University Magdeburg\\
D-39016 Magdeburg\\ Germany\\
}

\maketitle

\begin{abstract}
A new almost perfect nonlinear function (APN) on $\F{10}$ which is 
not equivalent to any of the previously known 
APN mappings is constructed.  This is the first example of an APN
mapping 
which is 
not equivalent to a power mapping. \\ 


\end{abstract}

\section{Introduction}
In cryptography, one is interested in functions $F:\F{m}\to\F{m}$
which are highly nonlinear. There are basically two concepts to measure
the linearity of a function: We may use the Walsh transform
(which is a special case of the Discrete Fourier transform)
or we may use differential properties of $F$. These two 
concepts yield to the notion of almost bent (AB) and
 almost perfect nonlinear (APN)
functions. Not many examples of such functions are known, and 
it was an open problem to decide whether the list of
known APN and AB functions is complete. Moreover, all the examples
constructed so far have been equivalent to power mappings. In this paper
we discuss the mapping 
$$
F:\F{10}\to\F{10},\quad x\mapsto x^3+ux^{36}
$$
where $u$ is a suitable element  in the multiplicative
group $\Fm{10}$ of $\F{10}$, see Theorem \ref{thm1}. It turns out that these 
mappings are inequivalent to any power mappings, hence
they are new. This is the first example of a new APN
mapping for several years, see \cite{dobbertin-01}, and it is the first example
of a mapping which is inequivalent to any power mapping. Moreover,
the mapping is crooked or, in other words, differentially affine. 

We emphasize that our function is inequivalent to a power
mapping in the general way described in \cite{carlet-charpin-zinoviev-98} and \cite{buda-carlet-pott-wcc}.
It seems that not much attention has been
paid so far to the  question whether the known classes
of APN or AB  functions are
inequivalent in this general sense or not. We are not aware of
any reference that shows whether the known classes are 
equivalent or not. 

In this paper, we use a dimension argument (that has not 
been used before in order to distingish APN or AB functions)
to prove that the function mentioned above is 
really new. This argument is motivated by the dimension arguments
that are used in order to distinguish difference sets, and it
may be applied also to distinguish the known classes
of APN and AB mappings.

Throughout this paper, let $\F{m}$ denote the finite
field with $2^m$ elements. This field is also a vector
space $\Fv{m}$ of dimension $m$ over $\F{}$, or
simply an elementary abelian group of order $2^m$. 
The differential and the linear properties of 
a function $F$ are only related to the additive
structure of $\F{m}$ and have nothing to do with the
multiplicative group. However, in order to construct functions
with good linear and differential properties, we will use the multiplication in $\F{m}$.
For a description of the differential and linear properties of
functions $F$, it is enough to consider $F$ to
be a mapping between two abelian groups, no matter whether these
are the additive groups of finite fields or not.

The paper is organized as follows. In the next Section, we
describe the notion of AB and APN and crooked mappings.
In Section \ref{eq}, we discuss the problem to determine the
equivalence classes of functions. In the final Section, 
we apply the results of Section \ref{eq} to show that 
our new APN function is inequivalent to the known ones.
We conclude the paper with some interesting open problems and
related questions.

\section{Nonlinear functions}

Let $U$ and $V$ be arbitrary groups. 
If $F:U\to V$ is a function, then we define the {\bf graph} 
$G_F$ of $F$ as follows:
$$
G_F:=\{(x,F(x))\ :\ x\in U\}\subseteq U\times V.
$$
We define
$$
\delta_F(a,b):=|\{(x-y,F(x)-F(y))=(a,b)\ :\ x,y\in U\}|.
$$
Note that $\delta_F(a,b)$ is the number of solutions $(x,y)$  to the
equations
\begin{equation}\label{eqn}
x-y =a,\qquad
F(x)-F(y)  =  b
\end{equation}
or the number of solutions
$$
F(y+a)-F(y)=b.
$$
If $F$ is linear (hence $U$ and $V$ are the additive groups of vector spaces) then 
$$
\delta_F(a,b)\in\{0,|U|\}.
$$
A function $F$ is {\bf differentially highly nonlinear} if 
$$
{\cal D}(F):=\max_{a\in U, b\in V, (a,b)\ne (0,0)} \delta_F(a,b)
$$
is small. 

We are now going to describe the differential properties of 
$F$ in terms of group algebras. Let $G$ be an arbitrary 
multiplicatively written 
group, and let $\K[G]$ denote the group algebra of $G$
over the field $\K$. The group algebra consists of 
the formal sums
$$
\sum_{g\in G} a_g g,
$$
where $a_g\in\K$. We can define an addition
$$
\left(\sum_{g\in G} a_g g\right)+\left(\sum_{g\in G} b_g
g\right)=\sum_{g\in G} (a_g+b_g)g
$$
and a multiplication 
$$
\left(\sum_{g\in G} a_g g\right)\cdot\left(\sum_{g\in G} b_g
g\right)=\sum_{g\in G}\sum_{h\in G} (a_h\cdot b_{h^{-1}g})g
$$
In order to distinguish the addition in $\K[G]$ from the 
composition of elements in $G$, we prefer to write the
group multiplicatively if we use group algebra notation. 
However, the groups that we are really using
are 
always additively written. 

If $D\subseteq G$, we identify $D$ 
with the element $\sum_{g\in D} g$, which we denote, by abuse
of notation, $D$ again. Moreover, if $A=\sum_{g\in G} a_g g$, then
$A^{(-1)}:=\sum_{g\in G} a_g g^{-1}$. Using this notation, we obtain easily
the equation
$$
G_F G_F^{(-1)}=\sum_{(a,b)\in U \times V} \delta_F(a,b) (a,b).
$$
This  shows that the $\delta_F(a,b)$'s are
the coefficients of the elements in $G_FG_F^{(-1)}$. The set (or
multiset if
we also count multiplicities)
$$
\{\delta_F(a,b)\ : \ (a,b)\in U \times V\}
$$
is called the differential spectrum of $F$. 

Characters are an important concept in the theory of
group algebras. We restrict ourselves
to abelian groups, otherwise we have to replace characters by 
higher dimensional  representations. Characters are precisely the
one-dimensional representations of a group $G$.

A {\bf character} is a homomorphism $G\to\C^\ast$. In the abelian case,
there are $|G|$ characters which form a group $\hat{G}$ which is
isomorphic with $G$. The transformation 
$$
\C[G]\to \C^{|G|},\quad 
\sum a_g g \to (\sum a_g\chi(g))_{\chi\in \hat{G}}
$$
is called the {\bf discrete Fourier transform}. 
The character that maps every group element to $1$ is 
denoted $\chi_0$.

We look at the case of functions $F:U\to V$ and the
elements $G_F\in \C[U\times V]$. If $F$ is linear or affine
linear, then
$G_F$ is a coset of a subgroup of $U\times V$, and then 
$$
\chi(G_F)\in\{0,\pm |U|\}.
$$
This follows from the well known
orthogonality relations for characters.
Therefore, it is natural to call a mapping {\bf highly nonlinear} 
if
$$
{\cal LN}(F):=\max_{\chi\in\widehat{U\times V}, \chi\ne \chi_0} |\chi(G_F)|
$$
is small.
The set (or multiset) of character values
$$
\{ \chi(G_F)\ : \ \chi\in \widehat{U\times V}\}
$$
is called the Fourier spectrum of $F$. We may define the Fourier and
the differential spectrum also for arbitrary sets $A\subseteq G$  or 
arbitrary group algebra
elements $A\in\C[G]$.

Now let us look at the special case of 
elementary abelian $2$-groups $U$ and $V$.
We return to the  general case of abelian groups in
the next section, since, in our opinion, the term ``equivalence
of functions'' is best explained in this more general context.

If $U$ and $V$ are elementary abelian $2$-groups, then $\delta_F(a,b)$ is always 
even hence we have
\begin{equation}\label{APN}
{\cal D}(F)\geq 2:
\end{equation}
The numbers $\delta_F(a,b)$ are even since the two equations (\ref{eqn}) 
have always an even number of solutions, you may just change $x$ and $y$.
We say that a function is {\bf almost perfect nonlinear (APN)} if $|U|=|V|$ and we have
equality in (\ref{APN}).

Similarly, one can show
\begin{equation}\label{B}
{\cal LN}(F)\geq 2^{|U|/2}.
\end{equation}
This can be proved easily using some well known properties
of the discrete Fourier transform.
If $|U|=|V|=2^{m}$, we have the improvement 
\begin{equation}\label{AB}
{\cal LN}(F)\geq 2^{(m+1)/2},
\end{equation}
see \cite{nyberg-94}.
 Functions which satisfy (\ref{AB}) with equality are
called {\bf almost bent (AB)}, whereas functions which satisfy
(\ref{B}) are called bent. Sometimes the term bent is
reserved just for the case of functions with $|V|=2$. 
It is well known that AB functions (which can exist
only in the case $m$ odd) are APN. 

The development of the concept  of nonlinearity
does not make use of finite fields. However, in order 
to construct examples of APN and AB functions it is useful to
equip $U$ and $V$ with the structure of a finite
field. In this case, we can describe our mappings by
polynomials. The degree of this polynomial will play an important role
in the next section.

The best studied functions are the power mappings $x^d$.
So far, all known constructions of APN and AB functions are related to  power mappings.
It has been checked at least up to $m=15$ (see \cite{dobbertin-2}) that the following table
gives a complete list of power APN mappings on $\F{m}$:

\begin{center}
\small{Table 1\\

Known APN power functions on $\mathbb{F}_{2^m}$.}

\begin{tabular}{|c|c|c|c|}
\hline
  & \scriptsize{Exponents $d$} & \scriptsize{Conditions} & \scriptsize{Reference}\\
\hline
\hline
\scriptsize{Gold functions} & \scriptsize{$2^i+1$} & \scriptsize{$gcd(i,m)=1$, $1\le i \le\frac{m-1}{2}$} & \scriptsize{\cite{gold-67},\cite{nyberg-94}}\\

\hline
\scriptsize{Kasami functions} & \scriptsize{$2^{2i}-2^i+1$} & \scriptsize{$gcd(i,m)=1$, $1\le i \le\frac{m-1}{2}$}  & \scriptsize{\cite{kasami-71},\cite{janwa-wilson-93}}\\

\hline
\scriptsize{Welch function}  & \scriptsize{$2^t +3$ }&  \scriptsize{$m=2t+1$}  &\scriptsize{\cite{dobbertin-99welch}}\\
\hline
\scriptsize{Niho function}  &\scriptsize{$2^t+2^\frac{t}{2}-1$, $t$ even} & \scriptsize{$m=2t+1$} & \scriptsize{\cite{Dobbertin-99niho}}\\
 & \scriptsize{$2^t+2^\frac{3t+1}{2}-1$, $t$ odd} &  & \\
\hline
\scriptsize{Inverse function } &\scriptsize{$2^{2t}-1 $}& \scriptsize{$m=2t+1$}& \scriptsize{\cite{nyberg-94},\cite{beth-ding-93}}\\
\hline
\scriptsize{Dobbertin function}  & \scriptsize{$2^{4i}+2^{3i}+2^{2i}+2^{i}-1$} & \scriptsize{$m=5i$} & \scriptsize{\cite{dobbertin-01}}\\
\hline
\end{tabular}
\end{center}

It turns out that, in the odd dimension case, the Gold, Kasami,
Welch and Niho functions are AB. The condition $i\leq \frac{m-1}{2}$
in the Gold
and Kasami case is not really restrictive: It just means that the
functions with $i > \frac{m-1}{2}$ are affine equivalent to those with 
$i\leq \frac{m-1}{2}$. For a thorough discussion of the notion of
equivalence, we
refer the reader to Section \ref{eq}. 

If a function is APN and bijective, then the inverse is also APN. The
inverse functions are not included in the table above.

\section{Equivalence of functions}\label{eq}

Let $D=\sum_{g \in G} a_g g$ be an arbitrary element in the
group algebra $\C[G]$, and let ${\cal L}$ be an 
automorphism of $G$. We define
$$
{\cal L}(D):= \sum a_g {\cal L}(g).
$$
Obviously, the differential spectra of $D$  and ${\cal L}(D)$
are the same, and also the Fourier spectra of 
$D$ and ${\cal L}(D)$ are the same. For the statement about the
Fourier spectrum note that the mapping
$$
\chi': g \mapsto \chi({\cal L}(g))
$$
is a character if $\chi$ is a character.
More generally, for $g\in G$, the elements $Dg$ and ${\cal L}D$
have the same differential and Fourier spectra.

Therefore it is natural to call two group algebra elements $D_1$ and
$D_2$ equivalent if there is an automorphism ${\cal L}$ of $G$
and a group element $g\in G$ such that ${\cal L}(D_1) = D_2g$. 

Now we want to specialize this concept of equivalence to functions.
This has been first done in \cite{carlet-charpin-zinoviev-98}, Proposition 3, therefore
we will call the notion of equivalence of functions
that stems from the notion of equivalence of group algebra elements 
{\bf CCZ equivalence}.

$F:U\to V$ and the corresponding group algebra elements 
$G_F$. The problem is that ${\cal L}(G_F)$ is not necessarily
a group algebra element that correponds to a function $F':U\to V$.

We call two functions $F_1:U\to V$ and $F_2:U \to V$ 
{\bf CCZ equivalent} if there is an automorphism ${\cal L}$ of $U \times V$
such that ${\cal L}(G_{F_1})=G_{F_2}\cdot g$ for some element
$g\in U\times V$. 
This generalizes the concept of affine equivalence. The
original definition of affine equivalence is as follows:

Let $U$ and $V$ be elementary abelian groups of order $2^m$, i.e.
the additive groups of verctor spaces over $\F{}$. We say that
$F_i:U\to V$, $i=1,2$ are affine equivalent if there are
linear mappings  ${\cal L}_1$ and ${\cal L}_2$  on $\F{m}$ and
elements $a\in U$, $b\in V$ auch that 
$$
F_2(x)={\cal L}_2(F_1({\cal L}_1(x+a)))+b.
$$

\begin{satz}
Two functions $F_1$ and $F_2$ are affine equivalent 
if and only if they are CCZ equivalent via an automorphism
${\cal L}$ of $U\times V$ such that ${\cal L}(V)=V$.
\end{satz}

Given two functions $F_1$ and $F_2$ it is sometimes  easy to decide
whether they are affine equivalent. It turns out  that the
algebraic degree of a function $F$ is an affine invariant
(we refer the reader to \cite{carlet-charpin-zinoviev-98}
for the precise definition of algebraic degree; it is the largest
2-weight of
the exponents that occur in the polynomial representation of $F$). 
But it seems that the important question whether affinely inequivalent
 functions
are CCZ equivalent has not been investigated. In \cite{buda-carlet-pott-wcc},
several classes of functions are constructed which are 
CCZ equivalent to the Gold power mapping but not 
affinely equivalent to any power mapping. This shows that
CCZ equivalence is really a coarser equivalence relation than
affine equivalence. As far as we know there is no proof that none
of the
APN mappings described in Table 1 are CCZ equivalent.

If $m$ is odd, it is known that the Fourier spectra
of  the inverse function and  
the Dobbertin function are different, hence these
two functions are not CCZ equivalent.  
Moreover, there spectrum has more than three values.
However, the Gold, Kasami, Welch
and Niho functions all have the same $3$-valued Fourier spectrum, hence they 
can be distinguished from the inverse and the Dobbertin function,
but they cannot be distinguished between themselves using the Fourier
spectrum.

In the case $m$ even, the Fourier spectrum of the Gold and Kasami
power functions are equal. It is
always different from the spectrum of the Dobbertin function,
see \cite{canteaut-charpin-dobbertin-00}.
It turns out that our new function has the same 
spectrum as Gold and Kasami. 
Therefore, in order to decide whether APN functions are
CCZ equivalent or not,  we have
to find other invariants than just the
Fourier spectrum. 

If $F:U \to V$ is an APN mapping (i.e. $U$ and $V$ are 
elementary abelian groups of order $2^m$), we define
$$
A_F:=\frac{G_F^2-2^m}{2}\in \C[G].
$$
We have $(a,b)\in A_F$ if and only
if
$F(x+a)+F(x)=b$ has two solutions in $x$. If $F_1$ and $F_2$ are
CCZ equivalent, then $A_{F_1}$ and $A_{F_2}$ are obviously equivalent.
Now we view $A_F$ as an element in $\F{}[U\times V]$.
If $F_1$ and $F_2$ are CCZ equivalent, then 
$A_{F_1}$ and $A_{F_2}$ are also equivalent in $\F{}[U\times V]$.
Hence the dimension of the ideal generated by $A_F$ 
is invariant under CCZ equivalence. 

We note that the ideal generated by $A_F$ may be also
viewed as the $\F{m}$ span of the following 
matrix ${\bf A}$ of size $2^{2m}\times 2^{2m}$: We index the
rows and columns with elements from $U \times V$. We have
$$
{\bf A}_{(a,b),(u,v)}=\left\{\begin{array}{ll}
1 & \mbox{if } (a+u,b+v)\in A_F\\
0 & \mbox{otherwise.}\end{array}\right.
$$

\section{The new APN function}

\begin{theorem}\label{thm1} 

Let $\omega$ be an element of order $3$ in $\F{10}$. 
Let $\F{5}$ denote the subfield of order $32$ in  $\F{10}$. The mapping 
\begin{equation}\label{ex1}
\begin{array}{cccl}
F: & \F{10} & \to & \F{10}\\
 & x & \mapsto &  x^3+u\cdot x^{36}
\end{array}
\end{equation}
is an APN mapping if and only if
\begin{equation}\label{cond1}
u\in \{\omega\Fm{5}\} \cup \{\omega^2\Fm{5}\}
\end{equation}
This function is not CCZ equivalent to any power mapping.
\end{theorem}

It is possible to give a ``theoretical'' argument why
these functions have the APN property.  Since this 
argument is quite involved, and since it does not really give insight
why the function is APN, we skip it. The APN property
of the function can be easily checked by computer.
One can easily show that the $62$ examples in (\ref{ex1}) are 
affine equivalent: In (\ref{ex1}), replace $x$ by $ax$ and then divide
the
resulting equation by $a^3$ to obtain
$$
x\mapsto x^3+\frac{u a^{36}}{a^3}x^{36} = x^3+ua^{33}x^{36}.
$$
But $ua^{33}$ satisfies (\ref{cond1}) if and only $u$ satisfies 
this condition (note $2^{10}-1=3\cdot 11\cdot 31$).

The function has the interesting property to be crooked. This
means that the sets
$$
H_a:=\{F(x+a)+F(x)\ :\ x\in\F{10}\}
$$
are affine hyperplanes in $\Fv{10}$. We refer the reader to
\cite{kyureghyan-wcc, kyureghyan-crooked}
for
recent progress on the problem to classify crooked mappings.

We want to distinguish our mapping from the known APN's. 
Table 1 shows that the only known APN mappings on $\F{10}$ 
are (up to affine equivalence)
$$
x^3, x^9\mbox{\ Gold }, x^{57}\mbox{\ Kasami }, x^{339}\mbox{\ Dobbertin }.
$$

As mentioned above, the Fourier spectrum of our new function is
different from the Fourier spectrum of the Dobbertin function.
This shows that our function  is
inequivalent to $x^{339}$.

The function $F$ is quadratic, therefore one may suspect  that
our function is affine or CCZ equivalent to one of the
Gold power mappings. 
Since the Fourier spectrum of our function is the same as those
of the Kasami and Gold function, we cannot use it to distinguish the
functions. 

We computed the dimensions of the ideals $I_F$  generated by $A_F$ 
for the Gold power functions $x^3$ and $x^9$ as well as the
Kasami power mapping $x^{57}$. The following table summarizes our
results:
\begin{center}
Table 2\\

Dimensions of the ideals $I_F$ in $\F{}[\Fv{10}\times \Fv{10}]$

\begin{tabular}{|c|c|}
\hline
$F$ & dimension\\
\hline \hline
$x^3$ & 1804\\
\hline
$x^9$ & 1804\\
\hline 
$x^{57}$ & $5734$\\
\hline
Theorem \ref{thm1} & 1896\\
\hline

\end{tabular}
\end{center}

This shows that our function is new. We can show that the
power mappings $x^3$ and $x^9$ on $\F{10}$ are not affine 
equivalent, but according to Table 2, the dimensions of the corresponding ideals
are the same. This shows that the dimension can not always
be used as a criteria to distinguish mappings.

It has been checked (by computer) that in finite fields
of order $\F{m}$, $m\leq 15$, there are no more
power APN mappings besides those listed in Table 1. Therefore, our function is not CCZ equivalent to
any power mapping. It was known before (see \cite{buda-carlet-pott-wcc}) that 
there are functions which are not affine equivalent to 
any power mapping. Our example gives the first
APN mapping which is not CCZ equivalent to any power
mapping. 

We can also use another argument if we want to show just
affine inequivalence to power mappings different from the
Gold case: Our function is crooked, and it is
known that the only crooked power mappings are quadratic,
see \cite{kyureghyan-wcc}. Hence the only chance to be affine equivalent is 
equivalence to the Gold power mapping, since affine equivalence
preserves the property being crooked. This property is
not preserved by CCZ equivalence!

The example in Theorem \ref{thm1} has been found through a computer search for 
APN binomials  $x^{d_1}+ux^{d_2}$ on $\F{n}$. 
The search was complete in the range $n\leq 10$. Up
to affine equivalence, the example in Theorem \ref{thm1}
is the only new APN binomial. We also found an example
in $\F{12}$ where we can show that the function is
not affine equivalent to the Gold power mappings. 

\begin{theorem}\label{thm2}
The mapping 
$$
\begin{array}{cccl}
F': & \F{12} & \to &  \F{12}\\
 & x & \mapsto & x^3+u\cdot x^{528}
\end{array}
$$
is an APN mapping if and only if 
\begin{eqnarray*}
u & \in & \{x\in\F{12}: \mbox{order of $x$ is divisible by $45$ and divides
  $45\cdot 13$}\} \\
& & \qquad \cup\  \{x\in\F{12}: \mbox{order of $x$ is divisible by 7 and divides
  $3\cdot 7\cdot 13$}\}.
\end{eqnarray*}
\end{theorem}

The proof that the functions
are not affine equivalent to the Gold power mappings is 
rather involved and therefore omitted. We 
did not yet check the dimension of the
ideal generated by $A_{F'}$ since the ambient space
is too large (it has dimension $2^{24}$). We also found some more
examples of binomials in larger fields where we are not 
yet able to prove that they are 
affine inequivalent to the known APN functions.  

\section{Summary and open problems}

In this paper, we reported about two new examples of APN functions
in $\F{10}$ and $\F{12}$. Both examples are
quadratic, which implies that the functions are crooked. 
In both cases, the new functions are
not affine equivalent to any power mapping, and in one case
we know that the example is not
CCZ equivalent to the Gold power mappings. Using 
computer assistance, we computed the dimensions of
the ideals generated by $A_F$ for different functions $F$.
These dimensions showed that the function on
$\F{10}$ is different from all previously known
APN mappings. Since all APN power mappings on $\F{10}$ are known, 
our function is not CCZ equivalent to any
power mapping.

We want to finish with the following open problems:

\begin{itemize}
\item Show that the function in Theorem \ref{thm2}
is not CCZ equivalent to any of the known functions.

\item Try to generalize the examples. Perhaps, one
can also use sums of more than just two Gold power mappings.

\item Give a theoretical proof that our new functions are not
CCZ equivalent to the known ones.

\item Compute the ranks of the ideals generated by
$A_F$ or $D_F$ for the known classes of APN or AB mappings.

\item Show that the known APN or AB functions are not
CCZ equivalent.  

\item Find more invariants for CCZ equivalence.
\end{itemize}

\vspace*{1cm}
\begin{center}
{\bf  Acknowledgment}
\end{center}
We thank J\"urgen Bierbrauer for many helpful discussions.
He actually initiated  our search for new APN
mappings. 

\vspace*{1cm}
\begin{center}
{\bf  Note}
\end{center}
After finishing this paper and making it available as a preprint, 
a new infinite series of APN functions  has been constructed (L. Budaghyan,
C. Carlet, P. Felke and
G. Leander: An infinite class of quadratic APN functions which are not equivalent to power mappings,
{\bf http//eprint.iacr.org/2005/359}). The series
covers some
of the examples presented here.  

\def\cprime{$'$} \def\cprime{$'$} \def\cprime{$'$}
  \def\Dbar{\leavevmode\lower.6ex\hbox to 0pt{\hskip-.23ex \accent"16\hss}D}
  \def\Dbar{\leavevmode\lower.6ex\hbox to 0pt{\hskip-.23ex \accent"16\hss}D}

\end{document}